%

\magnification=\magstep1
\input amstex
\documentstyle{amsppt}
\topmatter
\title
Small forcings and Cohen reals
\endtitle
\author
Jind\v rich Zapletal
\endauthor
\affil
The Pennsylvania State University
\endaffil
\address
Department of Mathematics,
The Pennsylvania State University,
University Park, PA 16802
\endaddress
\email
zapletal\@math.psu.edu
\endemail
\abstract
We show that all posets of size $\aleph _1$ may have to add a Cohen real 
and develop some forcing machinery for obtaining this sort of results.
\endabstract

\endtopmatter

\document

\subhead
{0. Results}
\endsubhead

\proclaim
{Theorem A} Cons (ZFC) implies Cons (ZFC+every forcing of size $\aleph _1$ 
adds a Cohen real).
\endproclaim

The proof of the Theorem A is a routine iteration argument based on 
the following two theorems which the author considers interesting 
in their own right.
\proclaim
{Theorem B} If $P$ is an $\aleph _0$-distributive forcing of size 
$\aleph _1$ then there is a proper forcing $Q$ such that $Q\Vdash$
``$\check P$ is nowhere $\aleph _0$-distributive".
\endproclaim
\proclaim
{Remark 1} Actually we get $Q\Vdash$``RO($\check P$)=
RO($Coll(\omega ,\omega _1))".$ Also if CH holds then our 
$Q$ will be $\aleph _1$-centered.
\endproclaim
\proclaim
{Corollary 1} PFA proves that there are no $\aleph _0$-distributive 
forcings of size $\aleph _1.$
\endproclaim
\demo
{Proof} Assume for contradiction that PFA holds and $|P|=\aleph _1,$ 
$P$ is $\aleph _0$-distributive. The Theorem B yields a $Q$ 
proper, $Q\Vdash$``$\langle \dot O_n:n\in \omega \rangle$ is a 
family of dense open subsets of $P$ with empty intersection". 
For $n\in \omega,$ $p\in P$ we define $D_{n,p}=\{ q\in Q:
\exists p^\prime \leq p$ $q\Vdash$``$p^\prime \in \dot O_n"\}$ 
and $D_{\omega ,p}=\{ q\in Q:\exists n\in \omega$ $q\Vdash$
``$\check p\notin \dot O_n"\} .$ Now $\langle D_{n,p},D_{\omega ,p}:
n\in \omega ,p\in P\rangle$ is a family of $\aleph _1$-many open dense 
subsets of $Q$ and by PFA there is a filter $F\subset Q$ meeting all 
of them. The reader can verify that if $\dot O_n/F=\{ p\in P:\exists 
q\in F$ $q\Vdash$``$\check p\in \dot O_n"\}$ then $\langle \dot O_n/F:
n\in \omega \rangle$ is a family of open dense subsets of $P$ with 
empty intersection, contradiction.
\enddemo

\proclaim
{Theorem C} If $P$ is a forcing adding a real then there is a c.c.c. 
forcing $Q$ such that $Q\Vdash$ ``$\check P$ adds a Cohen real".
\endproclaim
\proclaim
{Remark 2} If $|P|=\kappa$ then we obtain $|Q|=\kappa$ and in any 
case our $Q$ satisfies the Knaster condition.
\endproclaim
\proclaim
{Corollary 2} If $P$ is nowhere $\aleph _0$-distributive then there 
is a c.c.c. $Q$ such that $Q\Vdash$``$\check P$ adds a Cohen real".
\endproclaim
\demo
{Proof} Let $P\Vdash$``$\langle \dot \alpha _n:n\in \omega \rangle 
\subset \kappa$ is a new $\omega$-sequence of ordinals". Then if 
$Q_0$ is the forcing adding $\kappa$ Cohen reals we have $Q_0\Vdash$
``$\check P$ adds a real", namely the real coding the sequence of 
$\dot \alpha _n$-th Cohen reals. By the Theorem C, there is a c.c.c. 
$\dot Q_1\in V^{Q_0}$ such that $Q_0\Vdash \dot Q_1\Vdash$``$\check P$ 
adds a Cohen real". $Q_0*\dot Q_1$ is c.c.c. and we are done.
\enddemo

\proclaim
{Corollary 3} $MA_{<\kappa}$ proves that all nowhere 
$\aleph _0$-distributive forcings of size $<\kappa$ add a Cohen real.
\endproclaim
This uses something specific from the proof of the Theorem C and 
so we postpone the demonstration to the Section 2.

\proclaim
{Question 1} What if we want to embed more complicated forcings? 
For instance, if $C_{\aleph _1}$ denotes the BA adding $\aleph _1$ 
Cohen reals, is it consistent that every forcing of size $\aleph _1$ 
embeds $C_{\aleph _1}?$
\endproclaim
To show that the mechanism of giving a positive answer to this question 
would be different from the proof of the Theorem A we give
\proclaim
{Example 1} $MA_{\aleph _1}$ is consistent with existence of a forcing 
of size $\aleph _1$ adding a real but not embedding $C_{\aleph _1}.$
\endproclaim
\proclaim
{Question 2} What if we want to embed Cohen real into bigger forcings? 
Our proof of the Theorem B is specific for $\aleph _1.$ Is it consistent 
that all forcings of size $\aleph _2$ add a Cohen real?
\endproclaim

\subhead
{1. Proof of Theorem B}
\endsubhead

We use two simple lemmas.
\proclaim
{Lemma 1} If $P$ is an $\aleph _0$-distributive forcing of size 
$\aleph _1$ then there is a normal tree $T$ of height $\omega _1$ 
such that $T\subset P$ is dense. (The order of $T$ is inherited from $P. )$
\endproclaim
\proclaim
{Lemma 2} Forcings of the form $R*\dot S,$ where $R$ is the Cohen 
real and $R\Vdash$``$\dot S$ is $\sigma$-closed", do not add new 
branches to trees of height $\omega _1$ in the ground model.
\endproclaim
\proclaim
{Remark 3} Lemma 1 is specific for $\aleph _1.$ Lemma 2 holds true for
 many other standard generic reals in place of $R.$ Generalizations are 
left to the reader.
\endproclaim
Granted the lemmas, we prove the Theorem B: fix $P,$ an 
$\aleph _0$-distributive forcing of size $\aleph _1.$ 
Due to Lemma 1, we can chose $T_0\subset P,$ a dense normal tree 
of height $\omega _1.$ We construct $Q=Q_0*\dot Q_1*\dot Q_2:$
\roster
\item [0] $Q_0$ is the Cohen real.
\item [1] $\dot Q_1\in V_{Q_0}$ is the following set: $\{ \langle s,
 A\rangle :$ $s\in [T_0]^{\aleph _0},$ $A$ is a countable set of 
cofinal branches of $T_0\} .$
\endroster
Ordering is defined by $\langle s_0, A_0\rangle \geq \langle s_1,
 A_1\rangle$ if $s_0\subset s_1,$ $A_0\subset A_1$ and $\forall b\in 
A_0$ $b\cap s_1\subset s_0.$ Obviously, $\dot Q_1$ is $\sigma$-closed 
in $V^{Q_0}.$ Now if $\dot G_0*\dot G_1$ is the canonical term for a 
$Q_0*\dot Q_1$-generic filter (with the obvious meaning) we set 
$T_1\in V^{Q_0*\dot Q_1}$ to be $\bigcup \{ s:$ $\exists q\in 
\dot G_1$ $q=\langle s, A\rangle$ for some $A\} .$
\proclaim
{Claim 1} $Q_0*\dot Q_1\Vdash$``$\dot T_1\subset \check T_0$ is a 
dense subtree without cofinal branches".
\endproclaim
Given the claim,
\roster
\item [2] $\dot Q_2\in V^{Q_0*\dot Q_1}$ is the standard c.c.c.
 forcing specializing $T_1$ (\cite {She}).
\endroster

Obviously, $Q\Vdash$``$P$ adds a cofinal branch through a special 
tree of height $\omega _1"$ and so $Q\Vdash$``$\check P$ collapses 
$\aleph _1".$ For future reference we record
\proclaim
{Lemma 3} If $Q$ is as above then
\roster
\item $Q$ is proper
\item (CH) $Q$ is $\aleph _1$-centered
\item (GCH) $|Q|=\aleph _2$
\item $Q$ has $\omega _2$-p.i.c. (\cite {She,Ch. VIII})
\endroster
\endproclaim

 As for the Remark 1, note that all complete BA's of density 
$\aleph _1$ collapsing $\aleph _1$ are isomorphic to 
RO($Coll(\omega ,\omega _1)).$ 
\demo
{Proof of the Lemma 1} Let $P=\langle p_\alpha :\alpha <\omega _1\rangle$ 
be an $\aleph _0$-distributive forcing. We construct $T\subset P$ by 
induction on its levels. The induction hypothesis at $\alpha <\omega _1$ 
is that we have constructed $T^\beta \subset P,\beta <\alpha$ so that:
\roster
\item $T^\beta$ is a tree of height $\beta$
\item $\beta ^\prime <\beta <\alpha$ implies that 
$T^{\beta ^\prime }$ constitutes precisely the first
 $\beta ^\prime$ levels of $T^\beta$
\item the levels of $T^\beta$'s are maximal antichains in $P,\leq$
\item $\beta +1<\alpha$ implies $\exists t\in T^{\beta +1}$ $t\leq p_\beta .$
\endroster
How do we proceed with the induction?
\roster
\item  $\alpha$ limit. Set $T^\alpha=\bigcup _{\beta <\alpha }
T^\beta$ and the induction hypotheses continue to hold.
\item  $\alpha =\beta +1.$ As $P$ is $\aleph _0$-distributive, 
$D=\{ p\in P:\forall \beta ^\prime <\beta$ $\exists t$ in the 
$\beta ^\prime$-th level of $T^\beta$ $p\leq t\}$ is open dense in 
$P.$ (Remember (3)!) Choose $A\subset D,$ a maximal antichain such 
that there is $t\in A,$ $t\leq p_\beta .$ Set $T^\alpha =T^\beta \cup A$
 and the induction hypotheses again continue to hold.
\endroster
Finally, set $T=\bigcup _{\alpha <\omega _1} T^\alpha .$ Checking 
the needed properties of $T$ is trivial and we leave it to the reader.
\enddemo
\demo
{Proof of the Lemma 2}Fix $\dot S,$ $R\Vdash$``$\dot S$ is $\sigma$-closed", 
and $T,$ a tree of height $\omega _1.$ Assume $\langle r_0,\dot s_0\rangle 
\in R*\dot S,$ $\dot b$ are such that $\langle r_0,\dot s_0
\rangle \Vdash$``$\dot b$ is a new cofinal branch through $\check T".$
\proclaim
{Claim 2} For all $r_1,\dot s_1,t$ such that $r_1\leq r_0,$ $r_1\Vdash 
\dot s_1\leq \dot s_0$ and $\langle r_1,\dot s_1\rangle \Vdash$
``$\check t\in\dot b"$ there is $\alpha <\omega _1,$ $r_2\leq r_1$ 
and $\langle \dot s^i,t^i:i\in \omega \rangle$ such that $t^i$'s 
are distinct elements of the $\alpha$-th level of $T,$ $r_2\Vdash _R$
``$\dot s^i\leq \dot s_1"$ and $\langle r_2,\dot s^i\rangle\Vdash$
``$\check t^i\in \dot b".$
\endproclaim
\demo
{Proof} Let $r_1,\dot s_1,t$ witness the failure of the claim. Then 
immediately $T_0=\{ x\in T: \exists  \langle r_2,\dot s_2\rangle 
\leq \langle r_1,\dot s_1\rangle$ $\langle r_2,\dot s_2\rangle \Vdash 
\check x \in \dot b\}$ is a tree of height $\omega _1$ and all levels
 countable. Also  $\langle r_1,\dot s_1\rangle \Vdash$``$\dot b\subset 
\check T_0$ is a new cofinal branch". Now $R$ does not add new cofinal
 branches to $T_0$ and in $V^R,$ $S$ does not add new branches to $T_0$
 either, since $S$ is $\sigma$-closed and $T_0$ is still an 
$\omega _1$-tree there. So  $\langle r_1,\dot s_1\rangle \Vdash$
``$\dot b\in V",$ a contradiction.
\enddemo
Now fix $\theta$ large regular and $\langle M_i:i<\omega \rangle ,$ 
a sequence of countable submodels of $H_\theta$ such that $M_i
\subset M_{i+1},\dot b, \langle r_0,\dot s_0\rangle$ in $M_0$ and 
let $\alpha _i=M_i\cap \omega _1.$ Enumerate intersections of 
$\alpha _i$-th levels of $T$ with $M_{i+1}$ by $\langle x_i^j:
j\in \omega \rangle .$ Fix $G\subset R$ generic and work in $V[G].$ 

We let $f:\omega \to \omega$ be unbounded with respect to functions 
in $V.$ By the Claim 2, there are $\langle r_1,\dot s_1\rangle \leq 
\langle r_0,\dot s_0\rangle ,$ $r_1\in G,$ and $j_0\in \omega ,$ 
$j_0>f(0)$ such that $\langle r_1,\dot s_1\rangle \Vdash$
``$\check x_0^{j_0}\in \dot b".$ By elementarity we can find this 
 $\langle r_1,\dot s_1\rangle $ in $M_1.$ Now using this argument 
repeatedly together with the Claim 2, by induction on $i\in \omega$
 we can build a sequence $\langle \langle r_i,\dot s_i\rangle ,j_i:
i\in \omega \rangle$ such that $r_i\in G,$ $j_i>f(i),$ $\langle r_i,
\dot s_i\rangle$ in $M_i,$ $\langle r_{i+1},\dot s_{i+1}\rangle \leq 
\langle r_i,\dot s_i\rangle ,$ $\langle r_{i+1},\dot s_{i+1}
\rangle \Vdash$``$\check x_i^{j_i}\in \dot b".$ Since $\dot S/G$ 
is $\sigma$-closed, the decreasing sequence of conditions $\langle
 \dot s_i/G:i\in \omega \rangle$ has a lower bound. Let $\dot s_\omega$
 be an $R$-name for it.

Back in $V,$ let $\alpha =sup_{i\in \omega }\alpha _i$ and find any $x$ 
in the $\alpha$-th level of $T$ and $\langle r_{\omega +1},$ 
$\dot s_{\omega +1}\rangle \leq \langle r_0,\dot s_\omega \rangle$
 with $\langle r_{\omega +1},\dot s_{\omega +1}\rangle \Vdash$
``$\check x\in \dot b".$ Then we define $g:\omega \to \omega$ by 
$g(i)=j$ if $x_i^j>x.$ This function should not be bounded by any 
function in $V,$ since it is forced to be greater than our $f\in V[G].$ 
However, $g$ is clearly in $V,$ and we have obtained a contradiction.
\enddemo
\demo
{Proof of the Claim 1} Work in $V^{Q_0}.$ The density of $\dot T_1
\subset \check T_0$ is clear since countably many branches cannot 
cover all of $T_0\restriction t$ for any $t\in T_0.$ Let $p\in Q_1,$ 
$\dot b$ be such that $p\Vdash _{Q_1}$``$\dot b\subset \dot T_1$ is a
 cofinal branch". We distinguish two cases.
\roster
\item $\exists p_0\leq p$ $\forall p_1, p_2 \leq p_0,$ $t_1, t_2\in T,$ 
if $p_1\Vdash$``$\check t_1\in \dot b"$ and $p_2\Vdash$``$\check t_2\in
 \dot b"$ then $t_1,t_2$ are compatible. Then $c=\{ t\in T:\exists t_1
\leq t, p_1\leq p_0$ $p_1\Vdash$``$\check t_1\in \dot b"\}$ is a cofinal 
branch through $T$ and $p_ 0\Vdash$``$\dot b\subset \check c".$ So 
if $p_0=\langle s, A\rangle ,$ we may set $p_1=\langle s, A\cup \{ c\}
 \rangle$ to obtain $p_1\leq p_0,$ $p_1\Vdash$``$\dot T_1\cap \check
 c\subset \check s"$ and therefore $p_1\Vdash$``$\dot b\subset \check s$ 
and thus $\dot b$ is not cofinal", contradiction.
\item Otherwise. Then setting $\dot c=\{ t\in T_0: \exists t_1\leq t, t_1
\in \dot b\}$ we have $p\Vdash$``$\dot c\notin V$ is a cofinal branch 
through $\check T_0"$ contradicting the Lemma 2 for iteration $Q_0*Q_1$ 
and the tree $T_0.$
\endroster
We have a contradiction in both cases and the Claim is proven.
\enddemo
\demo
{Proof of the Lemma 3}\roster
\item $Q$ is an iteration of three proper forcings.
\item If CH holds, the centeredness of $Q$ follows from 
some cardinal arithmetic. Let $\langle \tau _\alpha :\alpha <\omega _1
\rangle$ enumerate the $Q_0$-names for elements of $[T_0]^{\aleph _0}
\in V^{Q_0}.$ Let $F_{\sigma ,\alpha ,q}=\{ \langle \sigma ,\langle 
\tau _\alpha ,\eta \rangle ,q\rangle :$  $\sigma \Vdash _{Q_0}$
``$\langle \tau _\alpha ,\eta \rangle \in \dot Q_1",$  $ \langle 
\sigma ,\langle \tau _\alpha ,\eta \rangle \rangle 
\Vdash _{Q_0*\dot Q_1}$``$\check q\in \dot Q_2"\}$ for 
 $\sigma \in {^{<\omega }2},$ $\alpha <\omega _1,$ 
$dom(q)\in [T_0]^{<\omega },$ $q:dom(q)\to \omega .$
 Then there are $\aleph _1$-many  $F_{\sigma ,\alpha ,q}$'s, 
each  $F_{\sigma ,\alpha ,q}\subset Q$ is a centered system and 
$\bigcup _{\sigma ,\alpha ,q} F_{\sigma ,\alpha ,q}\subset Q$ is
 dense and the $\aleph _1$-centeredness of $Q$ follows.
\item $Q_0\Vdash$``all cofinal branches of $\check T_0$ are in $V"$ and by 
GCH there are only $\aleph _2$-many of them. Again by GCH, one can 
enumerate all $Q_0$-names for pairs $\langle s,A\rangle$ as in the
 definition of $Q_1$ by $\langle \tau _\alpha :\alpha <\omega _2\rangle .$ 
Then $D=\{ \langle \sigma , \tau _\alpha  ,q\rangle :$  $\sigma \in 
{^{<\omega }2},$ $\alpha <\omega _2,$ $dom(q)\in [T_0]^{<\omega },$ 
$q:dom(q)\to \omega ,$ $\sigma \Vdash _{Q_0}$``$ \tau _\alpha  \in 
\dot Q_1",$  $ \langle \sigma , \tau _\alpha  \rangle 
\Vdash _{Q_0*\dot Q_1}$``$\check q\in \dot Q_2"\} \subset Q$ is dense. 
By the cardinal arithmetic and $\aleph _2$-c.c. of $Q$ we have 
$|RO(Q)|=\aleph _1$ and since $Q\subset RO(Q)$ we are done. 
\item Let us repeat what this means. \endroster
\proclaim
{Definition} \cite {She,Ch.VIII,\S 2} $Q$ has $\omega _2$-p.i.c. 
if for some $\theta$ large regular, $\Delta \in H_\theta ,$ 
{\it for every} $i<j<\omega _2,$ $q, h,$ $M_i,M_j$ countable 
submodels of $\langle H_\theta ,\in ,\Delta \rangle$ with 
$i\in M_i,$ $j\in M_j,$ $Q\in M_i\cap M_j ,$ $M_i\cap i=M_j\cap j,$ 
$M_i\cap \omega _2\subset j,$ $q\in Q\cap M_i,$ $h:M_i\to M_j$ an 
isomorphism which is identity on $M_i\cap M_j$ {\it there is} 
$q^\prime \leq q,$ a master condition for $M_i$ such that 
$q^\prime \Vdash$``$h^{\prime \prime }(\check M_i\cap \dot G)
=\check M_j\cap \dot G)".$
\endproclaim

To prove this we show that $Q_0,$ $\dot Q_1,$ $\dot Q_2$ have
 $\omega _2$-p.i.c. in the respective models and by \cite 
{She,Ch.VIII,Lemma 2.3} we will be finished. Certainly $Q_0$ has 
$\omega _2$-p.i.c. since for any isomorphism $h$ as in the definition 
of p.i.c. $h\restriction Q_0=id.$ For $\dot Q_1$ work in $V^{Q_0}$ and 
fix $q,h,M_i,M_j$ as in the definition of p.i.c. with 
$T_0\in M_i\cap M_j .$ Choose $q=q_0=\langle s_0,A_0\rangle
 \geq q_1=\langle s_1,A_1\rangle \geq \dots \geq q_i=
\langle s_i,A_i\rangle \geq \dots ,i\in \omega ,$ a strongly 
generic sequence for $M_i.$ Set $q_\omega =\langle s_\omega ,
A_\omega \rangle ,$ where $s_\omega =\bigcup _{i\in \omega }s_i$ 
and  $A_\omega =\bigcup _{i\in \omega }A_i.$ Also set 
$q_{\omega +1}=\langle s_\omega ,A_\omega \cup h^{\prime \prime }
A_\omega \rangle .$ Then $q\geq q_\omega \geq q_{\omega +1}$ and we 
claim that $q_{\omega _1}$ is what we are looking for. Certainly it 
is a master condition for $M_i.$ It is enough to show 
$q_{\omega +1}\Vdash$``$\check h^{\prime \prime }(\check M_i\cap \dot G)
\subset (\check M_j\cap \dot G)".$ Now if $r=\langle s,A\rangle \in M_i
\cap Q_1$ then $\exists q_{\omega +2}\leq q_{\omega +1}$ 
$q_{\omega +2}\Vdash$``$\check r\in \dot G"$ iff $q_{\omega +1}\Vdash$
``$\check r\in \dot G"$ iff $r\geq q_{\omega _1}$ iff $\forall b\in A$ 
$b\cap s_\omega \subset s$ iff $\forall b\in A$ $h(b)\cap s_\omega 
\subset s$ iff $h(r)=\langle h(s)=s,h(A)=h^{\prime \prime }A\rangle 
\geq q_{\omega +1}.$ Here the first and second equivalences hold by the 
strong genericity of $q_{\omega +1}$ and the fourth is due to the fact 
that $s\subset s_\omega \subset M_i\cap M_j$ and for $b\in M_i,$ 
$b\subset T_0$ a cofinal branch we have $b\cap M_i\cap M_j=h(b)\cap M_i
\cap M_j$ as $h$ is an isomorphism identical on $M_i\cap M_j.$ We are
 finished for $\dot Q_1$ and the case of $\dot Q_2$ is easy again:
 $\omega _2$-p.i.c. follows from the c.c.c. of $Q_2$ and from $|Q_2|
=\aleph _1.$ (Any $h$ as in the Definition has to be identical on 
$M_i\cap Q_2.)$
\enddemo

\subhead
{2. Proof of the Theorem C}
\endsubhead
Fix $P,\dot r$ such that $P\Vdash$``$\dot r\in {^\omega 2}\setminus V".$
 We define $Q$ as the set of ordered pairs $\langle f,g\rangle$ satisfying 
the following conditions:
\roster
\item $\exists n\in \omega$ $dom(f)=n.$ $n$ is called the {\it height} 
of the condition.
\item $\forall i<n$ $f(i)=\langle I_i,W_i\rangle ,$ where $I_i,i<n$ 
are subsequent finite intervals of $\omega$ and $W_i\subset {^{I_i}2}.$
\item $dom(g)\in [P]^{<\omega}.$
\item $\forall p \in dom(g)$ $g(p)\in {^{\leq n}2},$ $p$ decides 
$\dot r\restriction \bigcup _{i<lth(g(p))}I_i$ and $\forall \sigma 
\in {^n2}$ $g(p)\subset \sigma$ implies $\exists p^\prime \leq p$ 
$\forall i<n$ $p^\prime \Vdash$``$\dot r\restriction \check I_i\in 
\check W_i$ iff $\sigma (i)=1".$
\endroster

The order is by coordinatewise extension.

{\it Motivation.} The first coordinates will generically compose a
 sequence $\langle I_i,W_i:i<\omega \rangle$ and the future Cohen 
real will then be read off $\dot r$ as $\dot c(i)=1$ iff 
$\dot r\restriction I_i\in W_i.$ The second coordinate is 
(approximately) a finite fragment of a future projection of
 $P$ into the Cohen real algebra.

\proclaim
{Lemma 4} $Q$ is c.c.c.
\endproclaim
\demo
{Proof} We aim for the Knaster condition of $Q.$ Let 
$\langle q_\alpha =\langle f_\alpha ,g_\alpha \rangle :
\alpha <\omega _1 \rangle$ be a sequence of conditions in $Q.$
 We can thin this sequence out to $\langle q_\alpha :\alpha \in 
S\rangle$ for some $S\subset \omega _1$ of full cardinality such 
that $|\{ f_\alpha :\alpha \in S\} |=1,$ $\{ dom(g_\alpha ):\alpha \in S\}$
 is a $\Delta$-system with root $r$ and $|\{ g_\alpha \restriction r:\alpha
 \in S\} |=1.$ (First we use countability of the set of candidates 
for $f_\alpha$'s, then a $\Delta$-system argument on $dom(g_\alpha )$'s
 and finally countability of the set of candidates for $q_\alpha 
\restriction r$'s.) By the definition of $Q,$ if $\alpha _0,$ 
$\alpha _1\in S$ then $q_{\alpha _0},q_{\alpha _1}$ are compatible: 
their common lower bound is $\langle f_{\alpha _0},g_{\alpha _0}
\cup g_{\alpha _1}\rangle =\langle f_{\alpha _1},g_{\alpha _0}
\cup g_{\alpha _1}\rangle .$ We are done.
\enddemo
For future reference notice that if $|P|=\aleph _1$ 
then $|Q|=\aleph _1$ and $Q$ has $\omega _2$-p.i.c.

Let $H\subset Q$ be generic. In $V[H],$ set $F=\bigcup
 \{ f:\exists \langle f,g\rangle \in H\} ,$ $G=\bigcup 
\{ g:\exists \langle f,g\rangle \in H\} .$

\proclaim
{Lemma 5}
\roster
\item $dom(G)\subset P$ is dense.
\item $\forall p_0,p_1\in dom(G)$ if $p_0,p_1$ are compatible in $P$ 
then $G(p_0),$ $G(p_1)$ are compatible elements of $^{<\omega }2.$
\item $\forall p_0\in dom(G)$ $\forall \sigma \in {^{<\omega }2}$
 if $G(p_0)\subset \sigma$ then $\exists p_1\leq p_0$ $p_1\in dom(G)$ 
and $\sigma \subset G(p_1).$
\endroster
\endproclaim

Granted the lemma, we show how in $V[H],$ $P$ adds a Cohen real: 
if $K\subset P$ is generic over $V[H],$ set $c=\bigcup _{p\in H
\cap dom(G)}G(p).$ (2) makes sure that this is a function. $c$ 
is Cohen over $V[H]:$ let $p\in P,$ $D\subset {^{<\omega }2}$ 
open dense. Using (1), find $p_0\leq p,$ $p_0\in dom(G).$ 
There is $\sigma \in D$ extending $G(p_0).$ By (3) we can find 
$p_1\leq p_0$ such that $\sigma \subset G(p_1);$ thus $p_1\Vdash _P$
``$\dot c$ meets $\check D"$ and by genericity we are done.

To verify the claim of the Corollary 3, note first that if $MA_{<\kappa}$ 
holds and $|P|=\lambda <\kappa$ and $P$ is nowhere $\aleph _0$-distributive 
then $P$ adds a real. (This is because $P$ adds a countable sequence 
to $\lambda$ and as $2^{\aleph _0}>\lambda$ this new sequence can be coded 
over $V$ by a real, which then has to  be new as well.) Now we know that 
then there is a c.c.c. poset $Q$ adding a function $G$ with properties 
described in the Lemma 4. It is a simple exercise in counting necessary 
open dense subsets of $Q$ to show that then $G$ with these properties 
exists in $V.$ The same proof as above then shows that in $V,$ $P$ 
adds a Cohen real.

\demo
{Proof of the Lemma 5}\roster
\item Let $q=\langle f,g\rangle \in Q,$ $dom(f)=n\in \omega$ and $p\in P.$ 
Find $p^\prime \leq p,$ $p^\prime$ deciding $\dot r\restriction
\bigcup _{i<n}I_i,$ where for $i<n$ $f(i)=\langle I_i, W_i\rangle .$
 If $\sigma \in {^n2}$ is such that $\sigma (i)=1$ iff $p^\prime \Vdash$
``$\dot r\restriction \check I_i\in \check W_i"$ one can easily check that 
$ q^\prime =\langle f,g\cup \{ \langle p^\prime ,\sigma \rangle \} \rangle 
\in Q,$ $q^\prime \leq q$ and $q^\prime \Vdash$``$p^\prime \in dom(G),"$ 
and density of $dom(G)$ follows.
\item Easy.
\item Choose $\langle f,g\rangle \in Q,$ $p\in dom(g)$ and $g(p)\subset 
\sigma \in {^{<\omega }2}.$ First we show
\endroster
\proclaim
{Claim 3} For every $n\in \omega$ there is $\langle f^\prime ,
g\rangle \leq \langle f,g\rangle$ such that $n\subset dom(f^\prime ).$
\endproclaim
\demo
{Proof} This can be proven by induction on $n\in \omega .$ Obviously 
for $n=0$ there is nothing to prove. Assume now that we have $\langle 
f^\prime ,g\rangle \leq \langle f,g\rangle$ with $n\subset dom(f^\prime ).$
 If $n+1\subset dom(f^\prime )$ then $\langle f^\prime ,g\rangle$ witnesses 
the claim for $n+1$ and the induction step follows. So assume 
 $n= dom(f^\prime ).$ For each $p^\prime \in dom(g),$ $\eta \in 
{^n2},$ $g(p^\prime )\subset \eta$ we choose $r_{\eta ,p^\prime}\leq 
p^\prime$ such that $r_{\eta ,p^\prime}\Vdash _P$``$\forall i<n$ 
$\sigma (i)=1$ iff $\dot r\restriction \check I_i\in \check W_i"$ where 
for $i<n$ $f^\prime (i)=\langle I_i,W_i\rangle .$ This is possible by 
the definition of $Q$ and the induction hypothesis. Now choose 
$$\bigotimes _{\tau \in 2}K_{\eta ,p^\prime ,\tau }\subset 
\bigotimes _{\tau \in 2}P_{\eta ,p^\prime ,\tau }$$
generic, where $P_{\eta ,p^\prime ,\tau }$ are just distinct copies 
of $P$ and $r_{\eta ,p^\prime}\in K_{\eta ,p^\prime ,\tau }.$ Our 
initial assumption about $\dot r$ now gives 
$$ \dot r/K_{\hat \eta ,\hat p^\prime ,\hat \tau}\notin 
V[\bigotimes _{\langle \eta ,p^\prime ,\tau\rangle \neq \langle 
\hat \eta ,\hat p^\prime ,\hat \tau \rangle }K_{\eta ,p^\prime ,\tau }]$$
for any $\hat \eta ,\hat p^\prime ,\hat \tau$ and thus one can find 
$I_n,$ a finite interval of $\omega$ starting at $\bigcup _{i<n}I_i,$ 
such that $\dot r/K_{ \eta , p^\prime ,\tau}\restriction I_n$ are all
 different elements of $^{I_n}2.$ This is possible since there are only 
finitely many reals to take care of. Here is the only place where we use 
the forced novelty of $\dot r.$ Now we set $W_n=\{ \dot r/K_{\eta ,
p^\prime ,1}:\eta \in ^n2, $ $p^\prime \in dom(g)$ and $g(p^\prime )
\subset \sigma \}$ and $f^{\prime \prime }=f^\prime \cup \langle n,
\langle I_n,W_n\rangle \rangle .$ The attentive reader can check that
 $\langle f^{\prime \prime },g\rangle \in Q$ and thus finish the induction 
step on his own.
\enddemo
Given the claim we can easily complete the proof of (3): let $n$ be the 
length of $\sigma .$ Choose $\langle f^\prime ,g\rangle \leq \langle 
f,g\rangle$ such that $n\subset dom(f^\prime )=m.$ Choose $\eta \in 
{^m2},$ $\sigma \subset \eta .$ By the definition of $Q$ there is
 $p^\prime \leq p,$ $p^\prime$ deciding $\dot r\restriction 
\bigcup _{i<n}I_i$ and such that $\eta (i)=1$ iff $p^\prime \Vdash$
``$\dot r\restriction \check I_i\in \check W_i"$ where for 
$i<n$ $f^\prime(i)=\langle I_i, W_i\rangle .$ Then as in (1) 
$\langle f,g\rangle \geq \langle f^\prime ,g\cup \{ \langle p^\prime ,
\eta \rangle \} \rangle \in Q,$ and since $p^\prime \leq p,$ $\sigma 
\subset \eta,$ (3) follows.
\enddemo

\subhead
{3. Proof of the Theorem A}
\endsubhead

Let us start with a model of GCH. Fix  $\langle x_\alpha :\alpha <
\omega _2 \rangle ,$ an enumeration with repetitions of objects of 
the form $x_\alpha =\langle A_z^\alpha :z\in \omega _1 \times \omega _1
\rangle ,$ where $A_z\in [\{ f:dom(f)\in [\omega _2]^{\aleph _0},$
 $rng(f)\subset \omega _2 \} ]^{\aleph _1}.$ By induction on 
$\alpha <\omega _2$ we build a countable support iteration 
$$ P=\langle P_\alpha :\alpha \leq \omega _2 ,\dot Q_\alpha :
\alpha <\omega _2 \rangle$$
together with sequences $\langle \tau _\alpha ^i:i<\omega _2\rangle$ 
with the following induction hypothesis: for $\beta <\alpha$
\roster
\item $P_\beta \Vdash$``$\dot Q_\beta$ is a proper $\omega _2$-p.i.c. 
poset of size $\aleph _2$ and we assume that the universe of 
$\dot Q_\beta$ is $\check \omega _2"$
\item $|P_\beta |=\aleph _2,$ $P_\beta$ is $\aleph _2$-c.c., 
$P_\beta \Vdash GCH$
\item $\langle \tau _\beta ^i:i<\omega _2\rangle$ is an 
enumeration of $P_\beta$-names for elements of $\dot Q_\beta$ 
(resp. elements of $\omega _2).$ Moreover, defining $\prec _\beta 
\subset \omega _1\times \omega _1$ in $V^{P_\beta }$ by $\gamma _0
\prec \gamma _1$ iff there is $f\in A_{\langle \gamma _0,\gamma _1
\rangle }^\beta$ such that $g$ given by $g(\delta )=
\tau _{f(\delta )}^\delta$ if $\delta \in dom(f)$ and $g(\delta )=1$ 
otherwise, we have
\item if $P_\beta \Vdash$``$\prec _\beta$ is an $\aleph _0$-distributive 
poset"  then $\dot Q_\alpha \in V^{P_\alpha }$ is any proper 
$\omega _2$-p.i.c. forcing of size $\aleph _2$ such that 
$P_\beta *\dot Q_\beta \Vdash$``RO($\langle \omega _1,
\prec _\beta$ )=RO($Coll(\omega ,\omega _1))"$ 
(see the Theorem B and Lemma 3).
\item  if $P_\beta \Vdash$``$\prec _\beta$ is a 
poset adding a real" then $\dot Q_\beta \in V^{P_\beta }$
 is any c.c.c. $\omega _2$-p.i.c. forcing of size $\aleph _1$ 
such that $P_\beta *\dot Q_\beta \Vdash$``$\langle \omega _1,
\prec _\beta \rangle$ adds a Cohen real" (see the Theorem C).
\item otherwise. Then $P_\beta \Vdash \dot Q_\beta =1.$
\endroster
For $\alpha$ limit \cite {She,Ch.VIII,\S 2} takes care about
 preservation of (2). By GCH and (2), there are only 
$\aleph _2$-many $P_\alpha$-names for elements of $\omega _2$ 
(or elements of $\dot Q_\alpha )$ and (3) continues to hold. 
For (1),(4),(5) there is nothing to check. The successor step
 is handled similarly.

Now by \cite {She} $P=P_{\omega _2}$ is a proper $\aleph _2$-c.c. 
notion of forcing. We show $P\Vdash$``all forcings of size $\aleph _1$
 add a Cohen real". Let $p\in P,$ $p\Vdash$``$x$ is a poset with
 universe $\omega _1".$ W. l. o. g. either
\roster
\item $p\Vdash$``$x$ is $\aleph _0$-distributive", or
\item $p\Vdash$``$x$ adds a real"
\endroster
since nowhere $\aleph _0$-distributive forcings af size 
$\aleph _1$ add reals. For the first case, by the $\omega _2$-c.c. 
of $P$ and preservation of $\aleph _1$ there is $dom(p)<\alpha <\kappa$ 
such that $p\Vdash_{P_\alpha }$``$\prec _\alpha$ is 
$\aleph _0$-distributive poset of size $\aleph _1"$ and 
$p\Vdash _P$``$x=\prec _\alpha /G\cap P_\alpha ".$ But then 
$p\Vdash _P$``$V^{P_{\alpha +1}}\models$ RO($x$)=
RO($Coll(\omega ,\omega _1))"$ and since the equality is 
absolute upwards as long as $\omega _1$ stays in place the 
same holds in $V^P.$ The second case is taken care of in the 
same way, observing that the formula ``$x$ adds a Cohen real" 
is absolute upwards.

This leaves us with only one thing to demonstrate, the Example 1. 
We define the following forcing $P:$ $P=\{ f: dom(f)\in \omega _1, 
rng(f)\subset {^{<\omega}2}\} .$ The ordering is defined by $f\geq g$ 
if $dom(f)\subset dom(g),$ $\forall \beta \in dom(f)$ $f(\beta )\subset 
g(\beta )$ and $\{ \beta \in dom(f):f(\beta )\neq g(\beta )\}$ is finite. 
As far as we know, $P$ has not been explicitly defined before, so we list 
some of its simplest properties:
\roster
\item The $P$-generic $G$ is unambiguously given by
 $F:\omega _1\to {^\omega 2},$ where $F(\beta )=\bigcup 
\{ f(\beta ):f\in G,\beta \in dom(f)\} .$ Each $F(\beta)$ 
is Cohen generic over the ground model.
\item $|P|=2^{\aleph _0}.$
\item $P$ is proper; actually, $P$ embeds into (Cohen subset 
of $\omega _1$ by countable conditions)$\times C_{\aleph _1}.$
\item $P$ embeds (Cohen subset of $\omega _1$ by countable conditions).
\item If $Q$ is c.c.c. then $Q\Vdash$``$C_{\aleph _1}$ does not 
embed into $\check P".$
\endroster

Only (2) and (5) are relevant for our purposes, and we leave the proof
of the other items to the reader. Notice that the consistency statement
 in Example 1 follows immediately: just start with $V\models$CH and 
force $MA_{\aleph _1}$ by a c.c.c. poset. (2) and (5) together show 
that in the resulting model $|P^V|=\aleph _1$ and $P^V$ does not 
embed $C_{\aleph _1}.$ Obviously, $P^V$ adds many new reals.

 Now (2) is trivial; we concentrate on proving (5). For contradiction, 
assume we have a c.c.c. forcing $Q,$ $q\in Q,$ $p\in P$ and $\dot h,$ 
a $Q$-name for a $P$-name such that $q\Vdash _Qp\Vdash _P$``$\dot h:
\omega _1\to 2$ is $C_{\aleph _1}$-generic over $V^Q".$ By induction 
on $\alpha <\omega _1$ we construct a sequence $\langle f_\alpha ,
s_\alpha ,t_\alpha ,i_\alpha ,q_\alpha :\alpha <\omega _1\rangle$ so that
\roster
\item $f_\alpha \in P,$ $s_\alpha \in [dom(f_\alpha )]^{<\omega},$ 
$t_\alpha :s_\alpha \to {^{<\omega }2},$ $i_\alpha \in 2$ and $q_\alpha 
\in Q,$ $q_\alpha \leq q.$
\item $f_0=p$ and the $f_\alpha$'s are continuously increasing with 
respect to ordinary inclusion. Also $\forall \beta \in s_\alpha$ 
$f_\alpha (\beta )\subset t_\alpha (\beta ).$ 
\item For two functions $k,l$ define $k\nearrow l$ to 
be $\{ \langle x,y\rangle :x\in dom(k)\setminus dom(l),$ $k(x)=y$ or
 $x\in dom(l),$ $l(x)=y\} .$ Then for each $\alpha <\omega _1$ we want
 $q_\alpha \Vdash _Qf_{\alpha +1} \nearrow t_\alpha\Vdash _P$
``$\dot h(\check \alpha )=i_\alpha ".$
\endroster
There is no problem in the induction. Once we are done, by a 
Fodor-style argument we find stationary $S\subset \omega _1$ 
such that $|\{ s_\alpha :\alpha \in S\} |=1,$ 
$|\{ t_\alpha :\alpha \in S\} |=1.$  Now $Q$ is c.c.c. and so there 
is $q^\prime \leq q,$ $q^\prime \Vdash _Q$``$|\{ \alpha \in S: 
q_\alpha \in \dot K\} |=\aleph _1",$ where $\dot K$ is the term 
for a $Q$-generic. Once more by c.c.c.-ness of $Q$ there is 
$\beta <\omega _1$ such that $q^\prime \Vdash _Q$``$\dot Z=\{ \alpha 
\in S\cap \beta :q_\alpha \in \dot K\}$ is infinite". Now set 
$p^\prime \in P,$ $p^\prime \leq p$ to be $f_\beta \nearrow t,$ where 
$t$ is the only element of $\{ t_\alpha :\alpha \in S\} .$ Then 
$q^\prime \Vdash _Qp^\prime \Vdash _P$``$\forall \alpha \in \dot Z$ 
$\dot h(\alpha )=i_\alpha "$ and so $q^\prime \Vdash _Qp^\prime 
\Vdash _P$``$\dot h\restriction \dot Z\in V^Q".$ Since $\dot Z\in V^Q$
 is an infinite set this contradicts our assumption about $\dot h$ 
being $C_{\aleph _1}$-generic over $V^Q.$

\enddocument